%

\magnification\magstep1
\hfuzz9pt

\input mssymb


%

%


\newcount\skewfactor
\def\mathunderaccent#1#2 {\let\theaccent#1\skewfactor#2
\mathpalette\putaccentunder}
\def\putaccentunder#1#2{\oalign{$#1#2$\crcr\hidewidth
\vbox to.2ex{\hbox{$#1\skew\skewfactor\theaccent{}$}\vss}\hidewidth}}


\def\eqdef{\buildrel \rm def \over =}
\def\rest{\mathord{\restriction}}
\def\ms{\medskip}

\def\phi{\varphi}

\def\su{\subseteq}
\def\a{\alpha}
\def\b{\beta}

\def\l{\lambda}
\def\k{\kappa}

\def\om{\omega}

\def\lng{\langle}
\def\rng{\rangle}
\def\ov{\overline}
\def\sm{\setminus}
\def\cont{{2^{\aleph_0}}}

\def\twoL#1{{}^{<#1}2}
\def\twoto#1{{}^{#1}2}
\def\omL#1{{}^{<#1}\om}
\def\omto#1{{}^{#1}\om}

\def\dom{{\rm dom}}

\def\ran{{\rm  ran}}

\def\id{{\rm id}}

\def\aut{{\rm Aut}}


\def\endproof#1{\hfill  
{\parfillskip0pt$\smiley_{\hbox{{#1}}}$\par\medbreak}}

\def\imply{\Rightarrow}
\def\iff{\Leftrightarrow}
\def\proof{\smallbreak\noindent{\sl Proof}: }

\def\conc{\concatenate}
\def\concatenate{\widehat{\;}}

\def\init{\triangleleft}



\newbox\noforkbox \newdimen\forklinewidth
\forklinewidth=0.3pt   

\setbox0\hbox{$\textstyle\bigcup$}
\setbox1\hbox to \wd0{\hfil\vrule width \forklinewidth depth \dp0
			height \ht0 \hfil}
\wd1=0 cm
\setbox\noforkbox\hbox{\box1\box0\relax}
\def\unionstick{\mathop{\copy\noforkbox}\limits}
\def\nonfork#1#2_#3{#1\unionstick_{\textstyle #3}#2}
\def\nonforkin#1#2_#3^#4{#1\unionstick_{\textstyle #3}^{\textstyle  
#4}#2}

\setbox0\hbox{$\textstyle\bigcup$}
\setbox1\hbox to \wd0{\hfil$\nmid$\hfil}
\setbox2\hbox to \wd0{\hfil\vrule height \ht0 depth \dp0 width
				\forklinewidth\hfil}
\wd1=0cm
\wd2=0cm
\newbox\doesforkbox
\setbox\doesforkbox\hbox{\box1\box0\relax}
\def\nunionstick{\mathop{\copy\doesforkbox}\limits}

\def\fork#1#2_#3{#1\nunionstick_{\textstyle #3}#2}
\def\forkin#1#2_#3^#4{#1\nunionstick_{\textstyle #3}^{\textstyle  
#4}#2}

\font\circle=lcircle10

\setbox0=\hbox{~~~~~}
\setbox1=\hbox to \wd0{\hfill$\scriptstyle\smile$\hfill} 
\setbox2=\hbox to \wd0{\hfill$\cdot\,\cdot$\hfill} 

\setbox3=\hbox to \wd0{\hfill\hskip4.8pt\circle i\hskip-4.8pt\hfill}  


\wd1=0cm
\wd2=0cm
\wd3=0cm
\wd4=0cm

\newbox\smilebox
\setbox\smilebox \hbox {\lower 0.4ex\box1
		 \raise 0.3ex\box2
		 \raise 0.5ex\box3
		\box4
		\box0{}}
\def\smiley{\leavevmode\copy\smilebox}

\headline={\tenrm  
\number\folio\hfill\jobname\hfill\number\day.\number\month.\number\year}


\outer\long\def\ignore#1\endignore{}

\newcount\itemno
\def\itm{\advance\itemno1 \item{(\number\itemno)}}
\def\ritm{\advance\itemno1 \item{)\number\itemno(}}
\def\startitm{\itemno=0 }
\def\aitm{\advance\itemno1 

\item{(\letter\itemno)}}

\newcount\secno
\newcount\theono

\catcode`@=11
\newwrite\mgfile

\openin\mgfile \jobname.mg
\ifeof\mgfile \message{No file \jobname.mg}
	\else\closein\mgfile\relax\input \jobname.mg\fi
\relax
\openout\mgfile=\jobname.mg

\newif\ifproofmode
\proofmodetrue            

\def\@nofirst#1{}

\def\neusection{\advance\secno by 1\relax \theono=0\relax}
\def\neuchap{\secno=0\relax\theono=0\relax}

\neuchap

\def\labelit#1{\global\advance\theono by 1%
             \global\edef#1{%
             \number\secno.\number\theono}%
             \write\mgfile{\@definition{#1}}%
}





\def\ppro#1#2:{%
\labelit{#1}%
\smallbreak\noindent%
\@markit{#1}%
{\bf\ignorespaces {#2}:}}





\def\@definition#1{\string\def\string#1{#1}
\expandafter\@nofirst\string\%
(\the\pageno)}

\def\@markit#1{
\ifproofmode\llap{{ \expandafter\@nofirst\string#1\ }}\fi%
{\bf #1\ }
}

\def\h@markit#1{
\ifproofmode\edef\nxt{\string#1\ }%
{\tenrm\beginL\nxt\endL}
\fi%
{{\bf\beginL #1\endL}}
}
 ^^L
\def\labelcomment#1{\write\mgfile{\expandafter
		\@nofirst\string\%---#1}}

\catcode`@=12

\def\refrence#1#2:{\write\mgfile{\def\noexpand#1{#2}}%
\areference{#2}}
\def\areference#1{\medskip\item{[#1]} \ignorespaces}

\newcount\referencescount
\def\numrefrence#1{\advance\referencescount1
\edef#1{\number\referencescount}%
\write\mgfile{\def\noexpand#1{#1}}%
\areference{#1}}

\def\numericalreferences{\let\refrence\numrefrence}

\newcount\scratchregister
\def\simplepro{\scratchregister\theono\advance\scratchregister1 

\edef\scratchmacro{\number\secno.\number\scratchregister}%
\expandafter\ppro\csname\scratchmacro\endcsname}


\font\teneuf=eufm10
\font\seveneuf=eufm7
\font\fiveeuf=eufm5
\newfam\euffam
\textfont\euffam=\teneuf
\scriptfont\euffam=\seveneuf
\scriptscriptfont\euffam=\fiveeuf
\def\frak{\fam\euffam\teneuf}

\def\on{\rest}

\def\D{{\frak D}}

\def\compl{\neg}
\def\B{{\cal B}}

\def\cut{\cap}

\def\union{\cup}

\def\noteq{\not=}
\def\G{\Gamma}

\def\ite#1 {\item{(#1)}}

\let\on\restriction
\let\cut\cap
\def\fct{{}^\omega\!\omega}
\def\rang{\ran}
\def\dom{{\rm dom}}

\def\aside#1\par{\par{\leftskip=0.6\hsize\rightskip=0pt plus 2cm\sl  
#1\par}}

\proofmodefalse
\def\(#1){\{#1\}}

\font\bigfont cmbx10 scaled \magstep2
\font\namefont cmbx10 scaled \magstep2
{{
\obeylines

\everypar={\hskip0cm plus 1 fil}
\bf
\parskip=0.4cm

{\bigfont INFINITE HOMOGENEOUS}
{\bigfont BIPARTITE GRAPHS }
{\bigfont WITH UNEQUAL SIDES}

\bigskip
{\rm February 1992}
\bigskip

\parskip0.1cm

{\namefont Martin Goldstern}
2. Mathematisches Institut
FU Berlin, Arnimallee 3
1000 Berlin 33, Germany
{\tt goldstrn@math.fu-berlin.de}

\bigskip

{\namefont Rami Grossberg}
Department of Mathematics,
 Carnegie Mellon University
Pittsburgh, PA 15213, USA
{\tt  rg2g+@Andrew.cmu.edu}

\bigskip

{\namefont Menachem Kojman}
Institute of Mathematics 

Hebrew  University of Jerusalem, Givat Ram
91904 Jerusalem, Israel
{\tt kojman@math.huji.ac.il}

}
\vfill
\everypar{}
\rm
\leftskip2cm\rightskip2cm

ABSTRACT.  We call a bipartite graph {\it homogeneous} if every  
finite
partial automorphism which respects left and right can be extended to
a total automorphism.

 A $(\kappa,{\lambda} )$ bipartite graph is a bipartite graph with
left side of size $\kappa$ and right side of size ${\lambda}$. We  
show
that there is  a homogeneous $({\aleph_0},2^{\aleph_0} )$
bipartite graph  of girth 4
(thus answering negatively a question by Kupitz and
Perles), and that depending on the underlying set theory all
homogeneous $({\aleph_0}, \aleph_1)$ bipartite graphs may be
isomorphic, or there may be $2^{\aleph_1}$ many isomorphism types of
$(\aleph_0,\aleph_1)$ homogeneous graphs.

\footline{\hfill}
\global\pageno0

\eject
}

\headline{\tensl\hfill Goldstern, Grossberg, Kojman: Bipartite  
Graphs\hfill}

{\bf \S\number\secno. Introduction}

A homogeneous graph is one in which
every finite partial automorphism extends to a total automorphism.  
All
countable homogeneous graphs were classified in [\LachlanWoodrow],  
and
countable tournaments were classified in [\Lachlan]   (see also
[\Cherlin]). 

When looking at countable homogeneous {\it bipartite} graphs, one  
sees
that there are only five types of such graphs: complete
bipartite graphs, empty bipartite graphs, perfect matchings,
complements of perfect matchings and the countable random bipartite
graph.

In this paper we study the structure of uncountable homogeneous
bipartite graphs which have two sides of {\it unequal} cardinalities.  
We
must make the following demand on the notion of automorphism to admit
this class of graphs: a bipartite graph has a left and a right side,
and automorphisms preserve sides (this is necessary, as otherwise a
partial finite automorphisms which switches two vertices from the
different sides cannot be extended to a total automorphism).

We call
a bipartite homogeneous graph with a left side of cardinality $\k$  
and a
right side of cardinality $\l>\k$ and which is neither complete nor
empty, a $(\k,\l)$ saS graph. The name
should mean ``symmetric asymmetric'', where the symmetry is local,  
and
the asymmetry is global, in having a bigger right hand side. (The
demand that saS graphs are neither complete nor empty is to avoid
trivial cases).

\medskip

The paper is organized as follows: In Section 1 we classify
homogeneous bipartite graphs, and  remark that
 there are only five types
of {\it countable} homogeneous bipartite graphs. Then we prove the  
existence
of $(\aleph_0,2^{\aleph_0})$ saS graphs. The existence of such graphs
answers negatively the following question by J.~Kupitz and
M.~A.~Perles: is it true that in every connected locally 3-symmetric
(see below)
bipartite graph of girth 4 which is not a complete bipartite
graphs both sides are of equal cardinality? (Kupitz and Perles proved
that the  answer is ``yes'' if
the graph is finite).

In the second section we count the number of non isomorphic
$(\aleph_0,\aleph_1)$ saS graphs under the assumption of the weak
continuum hypothesis. We  prove that the weak
continuum hypothesis (i.e., $2^{\aleph_0}<2^{\aleph_1}$, which is a
consequence of the continuum hypothesis) implies that there are
$2^{\aleph_1}$ pairwise non isomorphic $(\aleph_0,\aleph_1)$ saS
graphs.

In the third section we 

show that $\neg$CH + MA implies that there is only one
$(\aleph_0,\aleph_1)$ saS graph up to isomorphism. These results
together show that the number of isomorphism types of
$(\aleph_0,\aleph_1)$ saS graphs is independent of ZFC, the usual
axioms of Set Theory.

Our interest in homogeneous bipartite graphs started when M.~Perles
introduced to us the question of the existence of a locally symmetric
infinite bipartite graphs of girth 4
with sides of unequal cardinalities. (See \continuumTheorem\ below.) 

We are
grateful to him for this, and not less for his careful reading of the
paper and his helpful suggestions.

The notation we use is mostly standard, but we nevertheless specify  
it
here.

\bigskip

\ppro \prel NOTATION:
\startitm
\itm    A {\it bipartite graph} is a triple $\Gamma=\lng L,R,E \rng = 
\lng L^{\Gamma} , R^{\Gamma}, E^{\Gamma} \rng$
such
that $L\cap R=\emptyset$, $L$ and $R$ are non-empty and $E\su 
\{\(x,y): x\in L, y \in R\}$. 




$L\cup R$ is the set of {\it vertices} 

of $\Gamma$, $E$ is the set of {\it edges}.  Members of  $L$  and $R$
are called  left and right vertices, respectively. 


Abusing notation, we sometimes write $v\in \Gamma$, instead of $v\in
L\cup
R$. 

Abusing notation even more, we may write $L\times R$ for $\{\{x,y\}:
x\in L, y\in R\}$. 

$\Gamma=\lng L,R, E\rng$ is a {\it subgraph} of $\Gamma'=\lng
L',R',E'\rng$ if $L \subseteq L'$, $R \subseteq R'$, $E \subseteq  
E'$.
It is called an {\it induced} subgraph if in addition $E = E'\cut
L\times R$. 

\itm A bipartite graph $\Gamma=\lng L,R,E\rng$ is {\it complete} if
for all $x\in L$, $y\in R$ we have $\(x,y)\in E$ and is called {\it
empty} if $E=\emptyset$. If $\Gamma=\lng
L,R,E\rng$, the {\it complement graph} of $\G$, is the graph whose
edge set is $L\times R  \setminus E$
\itm    If $\Gamma$ is a bipartite graph and $v\in \Gamma$, the set
$\Gamma(v)=\{u:u\in \G,\;\(v,u)\in E \}$ is
called the {\it set of neighbors} of $v$.   ${\Gamma}$ is called a
{perfect matching} iff ${\Gamma}(u)$ is a singleton for every $u\in
{\Gamma}$. 

\itm A {\it square} in a graph $\Gamma$ is a quadruple of distinct
vertices, $v_1,\cdots,v_4$ such that
$\(v_1,v_4)\in E$ and  $\(v_i,v_{i+1})\in E$ for $1\le i\le 3$. 

\itm A {\it partial homomorphism} between two graphs ${\Gamma}_1$,
${\Gamma}_2$ is a partial map $f:{\Gamma}_1 \to {\Gamma}_2$ with the
property that for all $x,y\in \dom(f)$:  $\{x,y\}\in E_1$ iff
$\{f(x), f(y)\} \in E_2$         

\itm 

%



%
A {\it partial isomorphism} between bipartite
graphs ${\Gamma}$ and ${\Gamma} '$ 

is a 1-1  partial map from $L^{\Gamma} \cup R^{\Gamma} $ into
$L^{\Gamma'}\cup R^{{\Gamma} '}$  which preserves left and right 

(i.e.,  $f[L^{\Gamma} ]\su L^{\Gamma'} $, $f[R^{\Gamma} ]\su
R^{\Gamma'} $) 

and
preserves edges and non edges (i.e.,
$\(u,v)\in E^{\Gamma} $ iff $\{f(u),f(v)\}\in E^{{\Gamma} '}$
 for all  $u,v\in \dom f$). 

 Such a partial isomorphism $f$ is called a
(total) isomorphism if $f$ is a bijection between the vertices of
${\Gamma}$ and ${\Gamma} '$. 

\item{} $f$ is called a (partial) automorphism of ${\Gamma}$ 

 if $f$ is a (partial) isomorphism of ${\Gamma}$ to ${\Gamma}$  
itself.
(So we only consider (partial) automorphisms which respect left and
right sides.)
 $\aut (\Gamma)$ is
the group of all
automorphisms of $\Gamma$. 

\itm    A bipartite graph $\Gamma$ is {\it locally $n$-symmetric} if
there
is some $H\su \aut(\G)$ such that for every
$v\in \G$ and every  two $n$-tuples of neighbors of $v$,
$x_1,\cdots,x_n$ and $y_1,\cdots,y_n$,
 there 

is an automorphism $\phi\in H$ such that $\phi(v)=v$ and
$\phi(x_i)=y_i$ for all $1\le i\le n$.
 In such a case we say that
$H$ acts on $\Gamma$ in a locally $n$-symmetric manner. 

\itm A bipartite graph $\Gamma$ is {\it homogeneous} if every finite
partial automorphism can be extended to
an automorphism. If $H\su \aut (\Gamma)$ has the property that for  
every
finite partial automorphism $f$ of $\Gamma$ there is an automorphism
in $H$ which extends $f$, we say that $H$ {\it acts homogeneously} on
$\Gamma$.



\medskip
Kupitz and Perles proved

\ppro \kupitz Theorem: If ${\Gamma}$  is a finite, connected  
bipartite
graph of girth 4 

which is not complete, and is locally $3$-symmetric,
then $|L|=|R|$.

\medskip

We shall also need some standard set theoretic notation:  $\om$ is  
the
set of all natural numbers. We use the 

convention that $n=\{0,1,\ldots,n-1\}$, namely that a natural number
equals the set of all smaller natural numbers. By $\omto\om$ we  
denote
all functions from $\om$ to $\om$ and by $\omL\om$ we denote all
finite sequences from $\om$. ${}^n \omega $ is the set of all
sequences of natural numbers of length $n$, i.e., functions from $n$
into $\omega$. For $\eta\in {}^n \omega $, $i\in \omega$ we let
$\eta\conc i$ be the sequence $\eta $ extended by $i$, i.e.,
$\eta\cup \{\lng n,i \rng\}$.

 The relation $\eta\init \nu$ between the
sequences $\eta$ and $\nu$ denotes that $\eta$ is an initial segment
of $\nu$.   $Ord$ is the class of
ordinals.  An ordinal is equal  to the sets of all smaller ordinals,
$\alpha = \{{\beta}\in Ord: {\beta} < \alpha\}$.

By $f[A]$ we denote the
{\it range} of the function $f$ when restricted to the set $A$. An
{\it $n$-tuple $\ov x$} of a set $A$ is an ordered subset
$\{x(1),x(2),\ldots,x(n)\}\su A$ of size $n$. By $|A|$ we denote the
cardinality (finite or infinite) of the set $A$. By $\dom f$ we  
denote
the {\it domain} of a function $f$ and by $\ran f$ its {\it range}.
 If $A\su \om$,  the {\it complement} of $A$ (in $\om$)
is the set $\compl A \eqdef\om\sm A$.

The symbol $\forall ^\infty x\in A$ means ``for all but finitely many
$x$ in $A$''.

\neusection

\beginsection{\S1. What Homogeneous Bipartite Graphs Exist?}

Let us classify all bipartite homogeneous graphs.
 Suppose 

$\G=\lng L,R,E\rng$
is homogeneous. If both sides are of cardinality 1, there are only  
two
 possibilities. Suppose then that $x\not=y$ are on the
same side (say $L$). If $\G(x)=\G(y)$, by homogeneity $\G(x)=\G(z)$  
for
every $z\in L$, or, in other words, there is a set $B\su R$ such that
$B=\Gamma(x)$ for all $x\in L$. If $B$  and $R\sm B$ are proper  
subsets
of $R$, an easy violation of homogeneity follows. Therefore
$\Gamma$ is either a complete or an empty bipartite graph.

Thus, if $\G$ is neither complete nor empty, it must be that $x=y\iff
\G(x)=\G(y)$ for every $x,y\in L$ and for every $x,y\in R$ (a graph
which satisfies this equivalence is called {\it extensional}).

 Let us
first assume that for some $x\in L$, $\G(x)$ is a finite subset of  
$R$
of cardinality $n$. By homogeneity, $\{\Gamma(x):x\in L\}=\{u:u\su
R,\;|u|=n\}$. If $n>1$ and $|R|>n+1$, this leads to a contradiction
(try mapping two $x$-s with $n-1$ common neighbors  to two
other $x$-s with $n-2$ common neighbors). If $|R|=n+1$, $\G$
is a complement of a perfect matching of size $2n+2$.  So we are left
with the case $n=1$. One possibility is that $R=\{u\}$, and in this
case $\Gamma$ is a complete bipartite graph. Otherwise, $\G$ must be
a perfect matching!  Similarly, if $\G(x)$ is co-finite for some  
$x\in
L$, then $\G$ is a complement of a perfect matching. All this applies
when $L$ is replaced by $R$.

We are left, then, with the case that every $x\in L$ has an infinite
co-infinite set of neighbors in $R$ and vise versa. In this case we
prove that $\G$ satisfies for every $k,l<\om$ the following property:
\ms
\item{$(*)_{k,l}$} For every distinct $x_0,\cdots,x_k,y_0,\cdots,y_l$
in $L$ (in 

$R$) there are infinitely many $u\in R$ (in $L$) such that $u\in
\Gamma(x_i)$ and $u\notin \Gamma(y_j)$ for $i\le k,j\le l$.

\proof
Given $x_0,\dots,x_k,y_0,\ldots,y_l\in L$, let us first prove that
there is at least {\it one} $u\in R$ which is a neighbor of every
$x_i$ and not a neighbor of every $y_j$ for $i\le k,j\le l$. Let
$v\in R$ be any vertex. Pick distinct $x'_0,\cdots,x'_k\in \Gamma(v)$  
and
$y'_0,\cdots,y'_l\notin \Gamma(v)$ from $L$.  This is possible, since
$\Gamma(v)$ is infinite co-infinite. Now find an automorphism $\phi$
that takes $x'_i,y'_j$ to $x_i,y_j$ respectively, and $u:=\phi(v)$ is
as we want. Next suppose that there are $x_i,y_j$ as above for which
there are only finitely many $u$ as above, and suppose, furthermore
that the number of such elements $u$ is minimal for this choice of
$x_i,y_j$. As $L$ is infinite, there is some $z\in L$, $z\not=x_i$  
and
$z\not=y_j$ for $i\le k$ and $j\le l$. Let $u$ be as above. If $u\in
\Gamma(z)$ let $z:=y_{l+1}$, and otherwise let $z:=x_{k+1}$ to obtain  
a
choice of $x_i,y_j$ with a smaller number of $u$, and a hence a
contradiction.

We call a bipartite graph which satisfies $(*)_{k,l}$ for all  
$k,l<\om$ {\it random}, and
state without proof:

\ppro \random Fact: Every countable random bipartite 

graph is homogeneous and has girth 4.

\ppro \unique Fact: Every two countable random
bipartite graphs are
isomorphic to each other.

\ms
For proofs see [\ErdosSpencer] p.98 or [\ChangKeisler] p. 93 and  
p.129.

\ppro\element Remark: As a consequence of \unique, the set of  
sentences
$\{(*)_{k,l}:k,l<\om\}$ is a set of axioms of a complete first-order
theory (see [\ChangKeisler] p.~113). Also,  the sentences in this
theory are exactly those sentences whose probability to hold in a
randomly chosen bipartite graph of size $2n$ tends to 1 when $n$  
tends
to infinity. 

\ms

Let us sum up what homogeneous bipartite graphs there are:

\startitm
\itm complete bipartite graphs and empty bipartite graphs.
\itm perfect matchings and complements of  perfect matchings.
\itm homogeneous random bipartite graphs.

Evidently, it is class (c) that deserves attention. By the remark
above, all members of class (c) are elementarily equivalent to each
other (i.e., satisfy the same first-order sentences).
We already mentioned that
the countable members of class (c) are all isomorphic to the  
countable
random bipartite graphs, so we might ask:

\ppro \quest Question: What uncountable homogeneous bipartite graphs
are there? As (a) and (b) are trivial, the question is what  
uncountable
members of class (c) are there?
\ms

We shall now show that there are homogeneous random graphs with
countable left side and uncountable right side. We call these graphs
$(\aleph_0,\k)$ saS graphs when the cardinality of their right side  
is
$\k>\aleph_0$. Recall that above we
showed that if a homogeneous bipartite graph is neither complete nor
empty then it is extensional. This implies in particular that $|L|\le
2^{|R|}$ and $|R|\le 2^{|L|}$. Therefore if in a homogeneous  
non-trivial
bipartite graph $|L|=\aleph_0$, we have an a priori bound of
$2^{\aleph_0}$ on $|R|$. We shall see that this bound is attained:

\ppro\continuumTheorem Theorem: There is an $({\aleph_0}, \cont)$ saS
graph.

\proof
The left side of our graph will be $\omega$, and the right side will
be a  set of functions in $\fct$. 

We will construct our graph as a projective 

limit, in some appropriate sense, of
a sequence $\lng \Gamma_n:n<\om\rng$ of  finite bipartite graphs.

We shall need the following notion:

\ppro\magicdef Definition: We say that $\Gamma'$ is a ``magic  
extension'' of
	$\Gamma$ if
\startitm
\itm $\Gamma$ is an induced subgraph of $\Gamma'$.
\itm Every finite partial automorphism of $\Gamma$ extends to a total
automorphism of $\Gamma'$.
\ms

E.~Hrushovski proved in [\Hrushovski] the  following theorem: 

\ppro\udi Theorem: Every finite graph has a finite magic extension.

Looking at the proof in [\Hrushovski] one can see that the same
	theorem is still true if we replace ``finite graph'' by
	``finite bipartite graph''.  Hence, we get the following  
fact:

\ppro \magic Fact: For every finite bipartite graph $\Gamma$ there is
a finite bipartite magic extension $\Gamma'$.

The fact follows also from a more recent result by Herwig on magic
extensions of finite relational structures see [\Herwig].
\ms
We remark here that it is only the finite case that needed a proof,
because it is standard and easy that every infinite (bipartite) graph
has a magic extension of the same cardinality. (see [\ChangKeisler,
p.~214ff])

\ms

\noindent {\it Proof of \continuumTheorem:} 

We define now the construction of the sequence $\lng
\Gamma_n:1\le n<\om\rng$. The graph $\Gamma_n=\lng L_n,R_n, E_n\rng$  
has a
left side $L_n$ which is an initial segment of $\om$ (a natural
number) and a right side $R_n\su \omto{n}$, a finite set of sequences
of natural numbers of length $n$. Let $L_1=\{0,1\}$ and $R_1=\{\lng
1\rng,\lng 2\rng\}$ and $E_1=\{(0,\lng1\rng),(1,\lng 2\rng)\}$.   
(This
will ensure that the graph we get at the end is neither empty nor
complete).

We demand:

\startitm
\itm $L_{2i+1}=L_{2i}$
\itm $R_{2i+1}=\{\eta\conc 1:\eta\in R_{2i}\}\union \{\eta\conc
2:\eta\in R_{2i}\}$ and for every $x\in L_{2i+1}=L_{2i}$ and $\nu\in
R_{2i+1}$, $\(x,\nu)\in E_{2i+1}\iff \({x,\nu\rest(2i)})\in
E_{2i}$

So at even stages we ``double'' the points of the right side. Put  
more
precisely,

we can define 

 $\rho_{2i}(\eta)= \eta\conc 1$ for all $\eta\in R_{2i}$,
$\pi_{2i}(\eta)=\eta\rest 2i$ for $\eta\in R_{2i+1}$, and we let  
$\pi_{2i}$ and
$\rho_{2i} $ be the identity on $L_{2i}$. 

Thus, 

 although $\Gamma_{2i}$ is not
an induced subgraph of $\Gamma_{2i+1}$, $\rho_{2i}$, is 

an embedding of $\Gamma_{2i}$ in $\Gamma_{2i+1}$ as an induced  
subgraph,
and $\pi_{2i}$ is a graph homomorphism.

At odd stages $2i+1$, we do the following:  Let $\rho_{2i+1}(\eta)=
\eta\conc 1$ for $\eta\in R_{2i+1}$, $\rho_{2i+1}= $ identity on
$L_{2i+1}$. Now find a magic extension  ${\Gamma}_{2i+2}= \lng
L_{2i+2}, R_{2i+2}, E_{2i+2}\rng$ of the graph
$\rho[{\Gamma}_{2i+1}]$. 

By renaming vertices we may assume that all vertices in $R_{2i+2} $
which   are not already in $\rho_{2i+1}[{\Gamma}_{2i+1}]$ are 

sequences of length $2i+2$ whose first $2i+1$ entries are all 0, and
that  $L_{2i+2}$ is an initial segment of the natural numbers.

Again we let $\pi_{2i+1}(\eta) = \eta\rest (2i+1)$ for all $\eta\in
\rho(R_{2i+1})$, $\pi(x)=x$ for $x\in L_{2i+1}$.  So $\pi$ is a  
partial
homomorphism from 

${\Gamma}_{2i+2}$ onto ${\Gamma}_{2i+1}$.




Note that our sequence of graphs, together with the maps $\pi_i$ can
be viewed almost as a projective system, except that the homomorphism
involved are partial.  Nevertheless, its ``projective limit'' can
be defined in a natural way:

We define ${\Gamma}_\infty =(L,R,E)$ as follows: 

 The left side $L=\omega$. The right side
$R=\{\eta\in\omto{\om}:(\forall^\infty n)(\eta\rest n\in R_n)\}$.
Let $E=\{\(x,\eta):(\forall^\infty n)(\(x,\eta\rest n)\in E_n)\}$.

We have to show two facts:

\ppro \card Fact: The cardinality of $R$ is $\cont$.
\ms
\ppro \hom Fact: The graph $\Gamma_\infty=\lng L,R,E\rng$ is  
homogeneous.

\ms
  The proof of the first being trivial, let us turn to the proof of
  the second. Suppose $f$ is a finite partial automorphism of
  $\Gamma$. We can find ${n_0}$ which is large enough such that 

  $(\dom f\union \ran f)\cap L \su L_{n_0}$ and such that for any
  $\eta_1\not=\eta_2$ in $\dom f\union \ran f$, $\eta_1\rest
  {n_0},\eta_2\rest {n_0}\in R_{n_0}$ and $\eta_1\rest
  {n_0}\not=\eta_2\rest {n_0}$, and such that for every $x,\eta\in
  \dom f\union \rang f$, $\(x,\eta)\in E\iff \(x,\eta\rest {n_0})\in
  E_{n_0}$. So for each $n\ge n_0$, $f$  induces a  (finite) partial
automorphism $f_{n}$ of 

  $\Gamma_{n}$: $f_{n}(\eta\rest n) = f(\eta)\rest n$ for all  
$\eta\in
  \dom(f)\cap R$, $f_n(x)=f(x)$ for $x\in\dom(f) \cap L$. 

Suppose without loss of
  generality that ${n_0}=2i_0+1$. Let $\bar f_{n_0} = f_{n_0}$. 

Now argue by induction on $n\ge n_0$
  to get a sequence of partial automorphisms $(\bar f_n: n\ge n_0)$
satisfying the following for all $n \ge n_0$:
\startitm
\itm $\bar f_n$ is a partial automorphism of ${\Gamma}_n$, and if $n>
n_0$, then $\bar f_n$ is total. 

\itm $\bar f_n$ extends $f_n$. 

\itm $\pi_n \circ \bar f_{n+1} = \bar f_n \circ \pi_n $

 Given
  ${\bar f}_{2i-1}$ ($i> n_0$), a partial automorphism on  
$\Gamma_{2i-1}$,
  we can find a 

  total automorphism $\bar f_{2i}$ of ${\Gamma}_{2i}$ extending
  ${\bar f}_{2i-1}$  (or more precisely, extending $\pi_{2i-1}^{-1}
  \circ f_{2i-1} \circ \pi_{2i-1}$). 

   Condition (2) will automatically be satisfied.

 Now we have to define ${\bar f}_{2i+1}$.  We must have ${\bar
     f}_{2i+1}\rest L_{2i+1}= {\bar f}_{2i}\rest L_{2i}$, so it
     remains to define ${\bar f}_{2i+1} \rest R_{2i+1}$.
To satisfy  condition (3), we require 

$$(*)\qquad\hbox{ if $\bar f_{2i}(x)=y$, then 

  ${\bar f}_{2i+1}[\{x\conc 1, x\conc 2\}] = \{y\conc 1, y\conc  
2\}$.}$$
     For $x$ in $\dom(f_{2i})\cap R_{2i}$, 

     exactly one of $x\conc 1$, $x\conc 2$
     is in $\dom(f_{2i+1})$ (by our assumption on $n_0$), so (2) and
     (3) uniquely determine  the  behaviour of $\bar f_{2i+1}$ on
     $x\conc 1$ and $x\conc 2$ in this  case. 

     For $\eta\notin \dom(f_{2i})$, we define $\bar f_{2i+1}(\eta) $
     arbitrarily      satisfying $(*)$.






 Having done the 

induction, let $F$ be defined of $\Gamma$ as follows: for $x\in \om$, 

$F(x)=y\iff (\forall^\infty n)({\bar f}_n(x)=y)$ and for $\eta\in R$, 

$F(\eta)=\nu\iff (\forall^\infty n)({\bar f}_n(\eta\rest n)=\nu\rest 
n)$.

We have to check that this indeed defines an automorphism. Note that
  all the ${\bar f}_i$ extend each other as far as the left side is
  concerned, and that whenever $\eta\in R_i$, $j< i$ and $\eta\rest j
     \in R_j$, then ${\bar f}_i(\eta)\rest j = {\bar f}_j(\eta\rest  
j)$. 

\relax From this property it is easy to see that all $F$ is
  well-defined on 

  the right side of ${\Gamma}$, and since all the ${\bar f}_i$ are
  automorphism, also ${\bar f}$ will be an
automorphism.\endproof{\card\hom\continuumTheorem}

We do mention one more thing: The proof actually gave us the  
following
property:
$$(**)\quad\vcenter{\hsize0.7\hsize\noindent
 for every finite partial automorphism $f$ of $\Gamma$ 

there is a locally finite automorphism $F$ of $\Gamma$ 

extending $f$.}$$
 By a locally finite automorphism we mean a permutation $F$ 

of $\Gamma$ with the property that for every finite $A\su \om$ there  
is a
finite $B\supseteq A$ such that $F\rest B\in Sym(B)$.

\ms

\ppro \generalizing Remark: (1) A similar proof shows the existence  
of
$(\kappa, 2^\kappa )$ saS graphs for any infinite cardinal $\kappa$. 

\hfil\break
(2) If $\kappa < {\lambda} ' \le {\lambda}$, and if ${\Gamma}$ is a
$(\kappa,{\lambda})$ saS graph, then it is easy to find an induced
subgraph ${\Gamma}'$ which is a $(\kappa, {\lambda} ')$ saS graph.




\neusection
\beginsection{\S2. The number of $(\aleph_0,\aleph_1)$ saS graphs
under weak CH}

In this section we handle the question of the number of the
isomorphism types of $(\aleph_0,\aleph_1)$ saS graphs. An obvious
upper bound is $2^{\aleph_1}$, the number of isomorphism types of
graphs of size $\aleph_1$.   In this section we show 

that if $2^{\aleph_0}<2^{\aleph_1}$, then this upper bound is
realized: there are $2^{\aleph_1}$   isomorphism types of
$(\aleph_0,\aleph_1)$ saS graphs. In the next section we show that if  
CH fails and MA
holds, then all $(\aleph_0,\aleph_1)$ saS graphs are isomorphic to
each other, namely there is a unique isomorphism type of
$(\aleph_0,\aleph_1)$ saS graphs.

\medskip
The idea of the first proof is as follows: we construct a family  
$\cal
G$ of $2^{\aleph_1}$
{\it different} saS graphs sharing  the same fixed countable 

left side. An isomorphism
between two saS graphs being determined by its action on the left
side, an isomorphism between two saS graphs in $\cal G$ is
really a permutation of the left side. There are $2^{\aleph_0}$
permutations of a given countable set, therefore there are at most
$2^{\aleph_0}$ members in every equivalence class of $\cal G$ modulo
isomorphism. Therefore it follows by $2^{\aleph_0}<2^{\aleph_1}$ that
there are $2^{\aleph_1}$ such classes.

The construction of many different saS graphs is done by iteratively
extending a countable random graph $\om_1$ many times, preserving
homogeneity and preserving the left side, in $2^{\aleph_1}$ many
different ways.

\medskip

\ppro \plusminus Notation: 

The left side of all graphs in this section will be $\omega$.  Since
we deal only with extensional graphs, we will identify a vertex in  
$R$
with its set of neighbors in $L$, so the edge relation will always be
given by $\in$. 

\hfil\break
 For $u\in R$ denote  $u^+\eqdef u$ and
$u^-\eqdef \neg u$.
\hfil\break
For a finite function $\sigma:R\to \{+,-\}$ we let
$B_\sigma=\bigcap_{u\in \dom \sigma}u^{\sigma(u)}$. If $\G$ is  
random,
then for every finite function $\sigma:R\to\{+,-\}$ the set  
$B_\sigma$
is infinite. 

\medskip

We now prove a few technical lemmas concerning the structure of the
automorphism group of a random bipartite graph, which will be used  
later in
extending countable random bipartite graphs:

\ppro \Alemma Lemma: Suppose that $\Gamma=\lng \om,R,\in\rng$ is  
random,
that $u_0,\ldots,u_k\in R$ and that $f,g\in\aut (\Gamma)$ are two
distinct automorphisms of $\Gamma$. Then there are $u,v\in R$, both
not in the list $u_0,\ldots,u_k$
 such 

that for every $x\in u\sm v$,
$f(x)\not=g(x)$.

\ms
What this lemma says is, that if two automorphisms are different,  
then
they are different on a definable  infinite
set of vertices: the set of all points which are connected to some  
$u$
and not connected to some $v$. Moreover, the $u$ and $v$ may be  
chosen
quite freely.

\proof
We may assume by applying $g^{-1}$ to $f$ and $g$, that $g=\id$. As
$f\not=\id$, there is some $x$ such that $f(x)\not=x$. As $\Gamma$ is
random, there are infinitely many $u\in \Gamma$ which satisfy $x\in  
u$ but
$f(x)\notin u$. Pick one such $u$ with the property that both $u$
and $f(u)$ are not in the list $u_0,\ldots,u_\k$ and set $v:=f(u)$.
For every $x\in u$, $f(x)\in v$. So if $x\in u\sm
v$, $f(x)\in v$, while $x\notin v$. In
particular, $f(x)\not=x$.\endproof\Alemma

\ppro \Acorollary Corollary: If $\Gamma$ is random,
$u_0,\ldots,u_k\in R$ and 

$g_1,g_2,\ldots,g_{l}\in \aut 
(\Gamma)$ then there is some finite function $\sigma:R\to\{+,-\}$  
such
that $\{u_0,\ldots, u_k\}\cap\dom \sigma=\emptyset$, and such that  
for
every $x\in B_\sigma$,
$g_1(x),g_2(x),\ldots,g_l(x)$ are 

$l$ distinct members of $\om$.

\proof Apply \Alemma\ iteratively ${l \choose 2}$
times.\endproof{\Acorollary}

\ppro \Blemma Lemma: Suppose that $B$ is an infinite subset of $\om$
and that $g_1,\ldots,g_k$ are 1-1 functions defined on $B$ with the
property that for every $x\in B$ and $1\le i < j\le k$, $g_i(x)\noteq
g_j(x)$.  Then there
is an infinite subset $B'\su B$ such that for every $x\not=y$ in
$B'$ and $1\le i\le j\le k$, $g_i(x)\not=g_j(y)$

\proof By induction on $n$ we pick an increasing chain of finite sets
$A_n$
with this property. At the induction stage: Clearly $g_i^{-1}[A_n]$  
is
finite, because $g_i$ is 1-1. Pick any $x\in B\sm \{g_i^{-
1}[g_j[A_n]]:1\le i\le j\le k\}$
and let $A_{n+1}=A_n\union \{x\}$. \endproof{\Blemma}

\ppro \something Corollary: Suppose that $\G$ is random,
$\tau:\G\to\{+,-\}$ is a finite partial function, and $G\su
\aut (\G)$ is finite. Then there exists an infinite set $B\su
B_\tau$ such that for $(g,x)\not=(g',x')\in G\times B$,
$g(x)\not=g'(x')$.

\proof By Corollary \Acorollary\ there is a finite function
$\sigma:R\to \{+,-\}$ such that for every $x\in B_\sigma$, the
elements $G(x)$ are $|G|$ distinct elements. As $\dom\tau$ and
$\dom\sigma$ are disjoint, also $\tau\cup\sigma$ is a function and
therefore $B_{\tau\cup\sigma}\su B_\sigma$ is infinite. By  Lemma
\Blemma \ there is an infinite set $B\su B_{\sigma\cup\tau}$ for  
which
the required holds.\endproof\something

We now prove the main lemma:

\ppro \extension Lemma:
Suppose $\Gamma$ is a countable random bipartite graph, and $G\su
\aut
(\Gamma)$ is a countable group of automorphisms. Then there are two
{\it incompatible} 

countable random bipartite graphs $\Gamma^0$ and $\Gamma^1$ with the
same left side as $\Gamma$, properly
extending $\Gamma$ such that $G\su
\aut (\Gamma^i)$ for $i\in \{0,1\}$.  By ``incompatible'' we mean  
that
there is no random
bipartite graph $\Gamma'$ with the same left side as $\Gamma$
extending both $\Gamma^0$ and $\Gamma^1$.

\proof
For any set $S\su\om$, the graph $\G^0=\lng \om,R\cup G(S),\in\rng$
satisfies that $G\su \aut (\G))$, where $G(S)$ denotes $\{g[S]:g\in
G\}$ and $g[S]=\{g(x):x\in S\}$. But $\G^0$ is not necessarily random
for an arbitrary choice of $S$ (for example, it is not random if
$S=\neg A$ for some $A\in R$).  We shall  find some subset $S$ of
$\om$ such that both $\G^0:=\lng \om,R\union G(S),\in\rng$ and
$\G^1=\lng \om,R\union G(\neg S)\in\rng$ are random.  This will
complete the proof, as in addition to the fact that $\G^i$ is random
and $G\su\aut \G^i$ ($i=0,1$), it is clear that there is no random
graph $\G'=\lng \om,R',\in\rng$ such that both $S$ and $\neg S$  
belong
to $R'$.

%

%

%

a
%

%

%

For $a,b \subseteq \omega $, $a\cut b=\emptyset$ and any finite
partial function ${\sigma} : G\cup R \to \{+,-\}$ we let 

$$
 B_{{\sigma},a,b}\eqdef 
  \bigcap_{A\in R\cut \dom({\sigma})}\!\! A^{{\sigma}(A)} \ \cap\ 
  \bigcap_{g\in G\cut \sigma^{-1}(+)}\!\! g(a)
  \bigcap_{g\in G\cut \sigma^{-1}(-)}\!\! g(b)
 $$
So our goal is to construct a set $S$ such that for all ${\sigma}$ as
above we have 

$$ B_{\sigma,S, \neg S} \hbox{ \ is  infinite }$$
This is equivalent to saying that  $\lng \om,R\union G(S),\in\rng$  
and by symmetry
also $\lng \om,R\union G(\neg S)\in\rng$ are random.

Note that there are only countably many finite partial functions
${\sigma}$ as above, so we can enumerate them as 

${\sigma}_1, {\sigma}_2, \ldots\,$.  We may also assume that each  
such
${\sigma}$ occurs infinitely many times in this list.

We construct $S$ and $\neg S$ by approximating them inductively
 by finite sets $a_n\su S$ and $b_n\su \neg S$, which satisfy
\item{(i)} $a_n\cap b_n=\emptyset$. 

\item{(ii)} $n\in a_{n+1}\cup b_{n+1}$
\item{(iii)} $a_n\su a_{n+1}$ and $b_n\su b_{n+1}$
\item{(iv)} There are at least $n$ points in  
$B_{\sigma_n,A_{n+1},b_{n+1}}$.

For $n=0$ let $a_n=b_n=\emptyset$.

For $n+1$: We specify which elements should be added to $a_n$ and
$b_n$ to obtain $a_{n+1}$ and $b_{n+1}$ respectively.

First, if $n\notin
a_{n}\cup b_{n}$, add it to $a_{n+1}$.



 We can find by \something\  an
infinite set $B\su B_{{\sigma}_n\rest R}$ such that for all
$(g,x)\not= (g',x')$ in  $(\dom({\sigma}\cut G))\times B$,
$g(x)\not=g'(x')$.

Now note that $ B_n\eqdef 
\{ g(x): x\in a_n\cup b_n, g\in \dom({\sigma}_n)\cut G \}$ is finite,
so we can find a set $X_n \su B \sm( B_n\cup\{n\}) $ of size $n$.

Let $a_{n+1}\eqdef a_n\cup \{g^{-1}(x): g\in {\sigma}^{-1}(+)\cut G,
x\in X_n\}$, and let $b_{n+1} \eqdef b_n\cup \{g^{-1}(x): g\in
{\sigma}^{-1}(-)\cut G, x\in X_n\}$. Note that $a_{n+1} $ and
$b_{n+1}$ are disjoint, by the conclusion of \something.  Moreover, $
B_{{\sigma}_n,a_{n+1} ,b_{n+1} } \supseteq X_n$ and therefore (iv) 

holds.

Let $S=\bigcup_n a_n$. It follows by (i) that $\neg S=\bigcup_n b_n$.
As for any $S \subseteq
\omega$, if $a_{n+1} \subseteq S$ and $b_{n+1} \subseteq S$, then $$
B_{{\sigma}_n,S, \neg S} \supseteq B_{{\sigma}_n,a_{n+1} ,b_{n+1} }
\supseteq X_n$$  Since each
finite partial $\sigma : G\cup R \to \{+,-\}$ appears as $\sigma_n$
with arbitrarily large $n$, this shows 

that $B_{{\sigma},S, \neg S}$ is infinite. 

\endproof{\extension}


%

%

$\sigma$ is a


\ppro \different Theorem: There are $2^{\aleph_1}$ different  
homogeneous
random bipartite graphs  of cardinality $\aleph_1$ with $\om$ as  
their
left side.

\proof
To  every $\eta\in \twoL{\om_1}$ we attach a pair
$\lng \Gamma_\eta,G_\eta\rng$ 

such that the following
conditions hold:

\startitm
\itm $\Gamma_\eta=\lng \om,R_\eta,\in\rng$ is a countable random
bipartite graph and $G_\eta\su \aut 
\Gamma_\eta$ is a countable group that acts on $\Gamma_\eta$  
homogeneously.
\itm If $\eta\init \nu$ then $R_\eta\su R_\nu$ and $G_\eta\su
G_\nu$.
\itm For every $\eta$,
$R_{\eta\conc 0}$ and $ R_{\eta\conc 1}$ are incompatible.

We define $\lng \Gamma_\eta$ and $G_\eta\rng$ by induction on
the length of $\eta$. 

If $\eta$ is the empty sequence, let
$\Gamma_\eta$ be any countable random bipartite graph with $\om$ as
its left side, and let $G_\eta$ be any countable group of
automorphisms that acts homogeneously on $\Gamma_\eta$.

If $\lg \eta$ is some limit ordinal $\a$, let
$R_\eta=\bigcup_{\b<\a}R_{\eta\rest\b}$ and let $G_\eta=\bigcup
_{\b<\a}G_{\eta\rest\b}$.  We should show that that $G_\eta\su
\aut(\Gamma_\eta)$ and that it acts homogeneously on $\Gamma_\eta$.  
As
all members of $G_\eta$ preserve $\in$ by their definition on  
$R_\eta$, it
is enough to show that $R_\eta$ is closed under $G_\eta$. Suppose  
that
$g\in G_\eta$ and $A\in R_\eta$ are arbitrary.  There is some
$\b<\lg(\eta)$ such that $g\in G_{\eta\rest \b}$ and $A\in
R_{\eta\rest\b}$. Now $g(A)\in R_{\eta\rest \b}\su R_\eta$.
 To see homogeneity, suppose $f$ is a finite partial
automorphism of $\Gamma_\eta$.  There is some ordinal $\b<\lg\eta$
such that $\dom f\union
\rang f\su \om\union
R_{\eta\rest \b}$. By the induction hypothesis, there is some $g\in
G_{\eta\rest \b}\su G_\eta$ extending $f$.

  If $\lng
\Gamma_\eta,G_\eta\rng$ is defined, use lemma \extension\ to find two
incompatible countable homogeneous random bipartite extensions of
$\Gamma_\eta$, $\Gamma_{\eta \conc 0}$ and $\Gamma_{\eta\conc1}$.  As
$\Gamma_{\eta\conc i}$ are {\it countable} random bipartite graphs  
for
$i\in \{0,1\}$, they are homogeneous by fact \random. For every  
finite
partial automorphism of $\Gamma_{\eta\conc i}$ there is an
automorphism of $\Gamma_{\eta\conc i}$ which extends it, so by adding
countably many automorphisms to $G_\eta$ and closing under  
composition
we get a countable group extending $G_\eta$ which acts
homogeneously on $\Gamma_{\eta\conc i}$. Let this group be
$G_{\eta\conc i}$.

Having done the definition by induction, we define for every sequence  
$\xi\in
\twoto {\om_1}$ a bipartite graph  
$\Gamma_\xi=\lng\om,\bigcup_{\a<\om_1}R_{\xi\rest \a},\in\rng$. 

As the group $G_\xi=\bigcup_{\a<\om_1}{G_{\xi\rest\a}}$ acts
homogeneously on $\Gamma_\xi$ --- as is easily seen --- $\Gamma_\xi$
is homogeneous.  Suppose that $\xi_0$ and $\xi_1$ are two different
members of $\twoto{\om_1}$ and let $\a$ be the last ordinal such that
$\xi_0\rest\a=\xi_1\rest \a$. By condition (3) above,
$\lng\om, R_{\xi_0},\in\rng$ $\lng \om,R_{\xi_1},\in\rng$ are
incompatible. As  $\Gamma_{\xi_0}$ and
$\Gamma_{\xi_1}$ are random, they must be different.
\endproof \different

\ppro\many Theorem: If $2^{\aleph_0}<2^{\aleph_1}$, then there are
$2^{\aleph_1}$ many isomorphism types of $(\aleph_0,\aleph_1)$ saS
graphs.

\proof
By the previous theorem there is a collection of $2^{\aleph_1}$ many
different saS graphs  $\{\Gamma_i:i<2^{\aleph_1}\}$ such that the
left side of each $\Gamma_i$ is $\om$. An isomorphism between  
$\Gamma_i$ and
$\Gamma_j$ for $i,j<2^{\aleph_1}$ is determined by its action on  
$\om$.
Therefore in an equivalence class of $\{\Gamma_i:i<2^{\aleph_1}\}$  
modulo
isomorphism there are at most $2^{\aleph_0}$ members. By the
assumption  $2^{\aleph_0}<2^{\aleph_1}$, it follows that there are
$2^{\aleph_1}$ many equivalence classes.\endproof \many

\ppro\kmany Remark: The proof above is readily generalized to give
$2^{\k^+}$ isomorphism types of $(\k,\k^+)$ saS graphs in case  
$2^\k<2^{\k^+}$.

We note that CH implies that $2^{\aleph_0}<2^{\aleph_1}$, and
therefore implies by the theorem above that there are $2^{\aleph_1}$  
many
isomorphism types of $(\aleph_0,\aleph_1)$ saS graphs.

\neusection
\beginsection{\S3. The number of  
$(\aleph_0,\aleph_1)$ saS graphs under the MA + $\neg$CH}

We
turn now to an examination of the number of $(\aleph_0,\aleph_1)$ saS
graphs under the assumption that CH fails but Martin's axiom MA  
holds. 

 The situation here is exactly opposite  to what we have seen under
(weak) CH. We shall prove the following:

\ppro\MAtheorem Theorem (MA): For any $\kappa< 2^{\aleph_0}$ there is
a unique  $({\aleph_0}, \kappa)$ saS graph.

First we will recall the statement of MA (see [\Kunen]):

A {\it dense set} $D$ in a partial order $(P,\le)$ is a subset $D\su
P$ such that for every $x\in P$ there is $y\in D$, $x\le y$.  Two
members $x,y\in P$ are {\it compatible} if there is $z\in P$ such  
that
$x\le z$ and $y\le z$. An {\it antichain} in $P$ is a set of pairwise
non compatible elements. A partial order satisfies the ccc (countable
chain condition) if every antichain is countable. A {\it filter} in a
partial order is a set $F\su P$ which satisfies (a) $F$ is {\it
downward closed}, i.e. $y\in F\;\&\;x\le y\imply x\in F$ and (b) $F$
is 

{\it directed}, i.e. $x,y\in F\imply (\exists z\in F)(z\ge  
x\;\&\;z\ge
y)$. The axiom MA (Martin's Axiom) is the statement ``for every ccc
partial order $P$ and every collection $\D$ of fewer than $\cont$
dense sets of $P$ there is a filter of $P$ with non-empty  
intersection
with every $D\in \D$''.  MA follows easily from the continuum
hypothesis (CH), but it is known
that MA is consistent with the negation of CH --- in fact, MA may be
true with the continuum being any regular cardinal.

\medskip

Let us introduce the following notation: if ${\Gamma} = (\omega,R,E)$  
is a
bipartite graph, ${\sigma}$ a finite partial function from $\omega$  
to
$\{+,-\}$ we let 

$$ \B_{\sigma}:=\{a\in R: \forall x\in
\dom({\sigma}): {\sigma}(x)={+}
\hbox{ iff } \(x,a)\in E\}$$

\ms

\ppro\bigsets Lemma: Let ${\Gamma}=\lng\om,R,E\rng$ be an
$(\aleph_0,\kappa)$ saS 

graph, $\kappa > {\aleph_0}$. Then for all ${\sigma}$ as above
we have $|\B_{\sigma}| = \kappa$.

\proof Fix $k$, $l$ in $\omega$.  We will only consider functions
${\sigma}$ with $|{\sigma}^{-1}(+)|=k$, $|{\sigma}^{-1}(-)|=l$.  For  
any
such functions ${\sigma}$, ${\sigma}'$ there is a partial  
automorphism
$f$ mapping $\sigma^{-1}(+)$ to ${\sigma'}^{-1}(+)$  and 

 $\sigma^{-1}(-)$ to ${\sigma'}^{-1}(-)$.  The total automorphism
$\bar f$ extending $f$ must map $\B_{\sigma}$ onto $\B_{\sigma} '$.

Hence all these sets $\B_{{\sigma}}$ have the same cardinality, say
${\lambda}$. Since 

since every element of $R$ must be in some such $\B_{\sigma}$ (by
homogeneity)  and there are only countably many such ${\sigma}$ we  
get
$\kappa \le {\lambda} \cdot {\aleph_0}$, i.e., ${\lambda} = \kappa$.
\endproof\bigsets

\ppro\getpartition Fact: If ${\Gamma}=\lng\omega ,R,E\rng$ is an
 $({\aleph_0},\kappa)$ saS graph, then $R$ can be partitioned into
$\kappa$ many countable 

sets  $(R_i:i<\kappa)$ such that for all $i<\kappa $ the
induced subgraph determined by $(\omega,R_i)$ is random.

\proof Let $R= \{x_i:i<\kappa\}$.  We will construct $(R_i:i<\kappa)$
by induction.  Given $(R_j:j<i)$, we can choose countable sets 

$$ R^{\sigma}_i \subseteq \B_{\sigma} \setminus \bigcup_{j<i} R_j$$
for every partial finite function ${\sigma}$ from $\omega$ to  
$\{+,-\}$,
because by \bigsets, $|\B_{\sigma}|=\kappa$, $ |\bigcup_{j<i} 
R_j|<\kappa$.  If $x_i\in  \bigcup_{j<i} R_j$ then let 
$$ R_i :=  \bigcup_{{\sigma}} R_j^{\sigma} $$
otherwise let $ R_i :=  \bigcup_{{\sigma}} R_j^{\sigma} \cup
\{x_i\}$. 

\endproof\getpartition

\ppro\forcingdef Definition: Assume   ${\Gamma}=\lng\omega ,R,E\rng$
and ${\Gamma} '= 
\lng \omega, R', E'\rng$  are  two $({\aleph_0},\kappa)$ saS graphs, 
and let $R=\bigcup_i R_i$, $R' = \bigcup_i R'_i$ be partitions as in
\getpartition. 

We let $P_{{\Gamma}, {\Gamma}'} $ be the set of all finite  partial
isomorphisms between ${\Gamma}$ and ${\Gamma} '$ respecting the
partitions, i.e., all finite partial isomorphisms $p$ satisfying 

$$ \forall x\in \dom(p)\cut R_i: p(x)\in R'_i$$
$P_{\Gamma,\Gamma'}$ is naturally ordered by the set inclusion
relation. (We consider functions to be  sets of ordered pairs.)

\ppro\cccFact Lemma: $(P_{{\Gamma}, {\Gamma} '}, {\subseteq })$ is a
partial order satisfying the countable chain condition.

\proof Let $\{p_\alpha:\alpha<\omega_1\} \subseteq P_{{\Gamma},
{\Gamma} '}$.  For each $\alpha$ let $s_\alpha:=\{i<\kappa:
\dom(p_\alpha)\cut R_i\not=\emptyset\}$.  $s_\alpha$ is a finite set. 

Applying the
${\Delta}$-system lemma [\Kunen, II, 1.5] we may without loss of
generality assume that $(s_\alpha:\alpha<\omega_1)$ forms a
${\Delta}$-system with root $s$.  Moreover, since there are only
countably many possibilities for $p_\alpha\on s$, we may also assume
that for some $p\in P_{{\Gamma}, {\Gamma} '}$ we have for all
$\alpha$: $p_\alpha\on s = p\on s$.  Similarly, we may assume
$p_\alpha\on \omega = p\on \omega$ for all $\alpha$.  Now for any  
$\alpha$,
${\beta}$ we have that $p_\alpha \cup p_{\beta}$ is a 1-1 function,
and hence an element of $P_{\Gamma,{\Gamma} '}$
\endproof\cccFact

\ppro\MAproof Proof of \MAtheorem: 

Let ${\Gamma}= \lng \omega, R,E\rng $,
${\Gamma} ' = \lng \omega, R', E'\rng $ be 

$(\omega, \kappa)$ saS graphs, and fix partitions as in  
\getpartition.
For any filter $G \subseteq P_{{\Gamma}, {\Gamma} '}$, we let $f_G:=
\bigcup G$.  Clearly $f_G$ will be a partial isomorphism from
${\Gamma}$ to ${\Gamma} '$.

Now note that for each $x\in \omega\cup R$, the set $D_x:= \{p\in
P_{{\Gamma}, {\Gamma} '}: x\in \dom(p)\}$  is a dense 

subset of $P_{{\Gamma}, {\Gamma} '}$ (because each $(\omega, R'_i)$  
is a
random bipartite graph).

By MA we can find a filter  $G \subseteq P_{{\Gamma}, {\Gamma} '}$
that meets all $D_x$.   This implies that $f_G$ is 

an isomorphism from ${\Gamma} $ into ${\Gamma} '$.  Similarly, using
$\kappa$ many dense sets defined from ${\Gamma}'$  we can insure that
$f$ will be onto.   Hence ${\Gamma}$ and ${\Gamma} ' $ are
isomorphic. \endproof\MAtheorem

conclude:


\vskip 2cm plus 0.5\vsize
\penalty-500
\vskip 0cm plus -0.5\vsize

\parindent1.5cm
\parskip\baselineskip
{\bf References:}

\def\paper#1,{{\sl #1},}

\refrence\Baumgartner B: J.~Baumgartner,\paper   All $\aleph_1$-dense  
sets of
reals can be isomorphic, {\bf Fundamenta Mathematicae} vol. 79  
(1971).
pp. 101--106.

\refrence\Cherlin C: G.~Cherlin, \paper  Homogeneous Tournaments  
revisited,
	{\bf Geometria Dedicata}  26 (1988), 231--240.

\refrence\ChangKeisler CK: C.~C.~Chang and H.~J.~Keisler, {\bf Model
	Theory}, North Holland, 1973.

\refrence\ErdosSpencer ES:
 P.~Erd\H os and J.~Spencer, {\bf Probabilistic Methods in
Combinatorics}, Academic Press 1974.

\refrence\Herwig HE: B.~Herwig, 

\paper Extending partial isomorphisms on finite structures, to
	appear in {\bf combinatorica}

\refrence\Hrushovski H: E.~Hrushovski, 

         \paper  Extending partial isomorphisms of graphs, to
	appear in {\bf combinatorica}

\refrence\Kunen K: K.~Kunen, {\bf Set Theory:  An
	introduction to independence proofs}, North Holland, 1980

\refrence\KupitzPerles KP:  Kupitz and Perles, \paper untitled, in  
preparation.

\refrence\Lachlan L:  A.~Lachlan, \paper  Countable Homogeneous  
Tournaments,
	{\bf Trans.\ Amer.\ Math.\ Soc.}\ 284 (1984), 431--461

\refrence\LachlanWoodrow LW: A.~H.~Lachlan and R~.E.~Woodrow,
	\paper  Countable Ultrahomogeneous Undirected Graphs, {\bf
Trans.\ Amer.\	Math.\ Soc.}\ 262 (1980) 51--94.

\bye